\documentstyle[12pt]{article}
\begin{document}

\author{S.V. Ludkovsky.}

\title{Operators on a non locally compact group algebra.}

\date{25 June 2012}
\maketitle

\begin{abstract}
The article is devoted to the investigation of operators on a non
locally compact group algebra. Their isomorphisms are also studied.
\footnote{key words and phrases: group, algebra,
operator, measure address: Department of Applied Mathematics, \\
Moscow State Technical University MIREA, \\ av. Vernadsky 78, Moscow
119454, Russia\\
Ludkowski@mirea.ru}

\end{abstract}

\section{Introduction.}
Group algebras play very important role in algebra, harmonic
analysis and operator theory
\cite{fell,fidal,hew,kawadamj48,losannm2008,nai}. Group algebras
were extensively studied for locally compact groups. One of the main
instruments in those investigations was an existence of a Haar
measure, which is characterized by such essential properties as of
being left or right invariant and quasi-invariant relative left and
right shifts and to the inversion on the entire group. \par But
substantially less is known for non locally compact groups. If a
nontrivial Borel measure on a topological Hausdorff group
quasi-invariant relative to the entire group is given, then such
group is locally compact according to A. Weil's theorem. Therefore,
on non locally compact Hausdorrf groups Borel measures may be
quasi-invariant relative to proper subgroups only. This is the
reason of many differences between group algebras of locally compact
and non locally compact groups. For non locally compact groups they
are already nonassociative. This work continues previous
publications of the author. \par In this article families of
topological groups which may be non locally compact are considered.
Group algebras of non locally compact Hausdorff topological groups
are studied. Particularly, operators on non locally compact group
algebras and their isomorphisms are investigated.  Borel regular
radonian measures $\mu _{\alpha }$ on topological groups $G_{\alpha
}$ quasi-invariant relative to dense subgroups $G_{\beta }$ are
taken. The Radon and Borel regularity properties for measures are
not very restrictive (see chapter 1 in \cite{dal} and chapter 2 in
\cite{federer}). The constructions of such measures were described
in \cite{dal,dalshn,lujms147:3:08,lujms150:4:08,luambp99,lunova2006}
and references therein. \par The main results of this paper are
obtained for the first time and are contained in Theorems 11, 15,
16, 18.

\section{Group algebra}
\par {\bf 1. Definition.} Let $\Lambda $ be a directed set and
$\{ G_{\alpha }: \alpha \in \Lambda \} $ be a family of topological
groups with completely regular (i.e. $T_1\cap T_{3\frac{1}{2}}$)
topologies $\tau _{\alpha }$ such that
\par $(1)$ $\theta ^{\beta }_{\alpha }: G_{\beta }\to G_{\alpha }$
is a continuous algebraic embedding with continuous inverse $(\theta
^{\beta }_{\alpha })^{-1}$, $\theta ^{\beta }_{\alpha }( G_{\beta
})$ is a proper subgroup in $G_{\alpha }$ for each $\alpha <\beta
\in \Lambda $;
\par $(2)$ $\tau
_{\alpha }\cap {\theta ^{\beta }_{\alpha }(G_{\beta })} \subset
\theta ^{\beta }_{\alpha }(\tau _{\beta })$ and $\theta ^{\beta
}_{\alpha } (G_{\beta })$ is dense in $G_{\alpha }$ for each $\alpha
<\beta \in \Lambda $;
\par $(3)$ $G_{\alpha } $ is complete relative to the left
uniformity with entourages of the diagonal of the form ${\cal U} =
\{ (h,g): h, g \in G_{\alpha }; h^{-1}g\in U \} $ with neighborhoods
$U$ of the unit element $e_{\alpha }$ in $G_{\alpha }$, $U\in \tau
_{\alpha }$, $e_{\alpha }\in U$;
\par $(4)$ for each $\beta =\phi (\alpha )$ the embedding $\theta ^{\beta }_{\alpha }$ is precompact, that is by our definition
for every open set $U$ in $G_{\beta }$ containing the unit element
$e_{\beta }$ a neighborhood $V\in \tau _{\beta }$ of $e_{\beta }$
exists so that $V\subset U$ and $\theta ^{\beta }_{\alpha }(V)$ is
precompact in $G_{\alpha }$, i.e. its closure $cl (\theta ^{\beta
}_{\alpha }(V))$ in $G_{\alpha }$ is compact, where $\phi : \Lambda
\to \Lambda $ is an increasing marked mapping.
\par {\bf 2. Definition.} Suppose that \par $(1)$ $\mu _{\alpha }: {\cal
B}(G_{\alpha }) \to [0,1]$ is a probability measure on the Borel
$\sigma $-algebra ${\cal B}(G_{\alpha })$ of a group $G_{\alpha }$
from \S 1 with $\mu _{\alpha }(G_{\alpha })=1$ so that \par $(2)$
$\mu _{\alpha }$ is quasi-invariant relative to the left and right
shifts on $h\in \theta ^{\beta }_{\alpha }(G_{\beta })$ for each
$\alpha <\beta \in \Lambda $, where $\rho ^{r}_{\mu _{\alpha
}}(h,g)=(\mu _{\alpha }^h)(dg)/\mu (dg)$ and $\rho ^{l}_{\mu
_{\alpha }}(h,g)=(\mu _{\alpha _h})(dg)/\mu (dg)$ denote
quasi-invariance $\mu _{\alpha }$-integrable factors, $\mu
^h_{\alpha }(S)=\mu (Sh^{-1})$ and $\mu _{\alpha ,h}(S)=\mu _{\alpha
}(h^{-1}S)$ for each Borel subset $S$ in $G_{\alpha }$. Moreover,
\par $(3)$ let a density $\psi _{\alpha }(g) =\mu _{\alpha
}(dg^{-1})/\mu _{\alpha }(dg)$ relative to the inversion exist and
let it be $\mu _{\alpha }$-integrable. \par A subset $E$ in
$G_{\alpha }$ has $\mu_{\alpha }$-measure zero, if a Borel subset
$F$ in $G_{\alpha }$ exists such that $E\subset F$ and $\mu _{\alpha
}(F)=0$. The completion of ${\cal B}(G_{\alpha })$ by all $\mu
_{\alpha }$-zero sets will be denoted by ${\cal A}(G_{\alpha })$.
The measure $\mu _{\alpha }$ has the extension $\nu _{\alpha }:
2^{G_{\alpha }}\to [0,1]$ such that $\nu _{\alpha } (E) := \inf \{
\mu _{\alpha }(F): ~ E\subset F \mbox{ and } F \in {\cal
B}(G_{\alpha }) \} ,$ where $2^{G_{\alpha }}$ denotes the family of
all subsets in $G_{\alpha }$. The measure $\nu _{\alpha }$ is Borel
regular, that is, by the definition all open subsets in $G_{\alpha
}$ are $\nu _{\alpha }$-measurable and each subset $E$ in $G_{\alpha
}$ is contained in a Borel subset $F$ so that $\nu _{\alpha }(E)=\nu
_{\alpha }(F)$. Evidently, $\nu _{\alpha }(F)=\mu _{\alpha }(F)$ for
each Borel subset $F$ in $G_{\alpha }$, so $\nu _{\alpha }$ on
$2^{G_{\alpha }}$ will also be denoted by $\mu _{\alpha }$.
\par Henceforth, it will be supposed that \par $(4)$ a subset $W_{\alpha }\in
{\cal A}(G_{\alpha })$ exists such that $\rho ^{r}_{\mu _{\alpha
}}(h,g)$ and $\rho ^{l}_{\mu _{\alpha }}(h,g)$ are continuous on
$\theta ^{\beta }_{\alpha }(G_{\beta })\times W_{\alpha }$ and $\psi
_{\alpha }(g)$ is continuous on $W_{\alpha }$ with $\mu _{\alpha
}(W_{\alpha })=1$ for each $\alpha \in \Lambda $ with $\beta = \phi
(\alpha )$. Let also each measure $\mu _{\alpha }$ be radonian, that
is for each $\epsilon >0$ a compact subset $V$ in $G_{\alpha }$
exists such that $\mu _{\alpha }(G_{\alpha }\setminus V)<\epsilon $.
\par {\bf 3. Notation.} Denote by $L^1_{G_{\beta }}(G_{\alpha },\mu_{\alpha })$
a complex subspace in $L^1(G_{\alpha }, \mu _{\alpha },{\bf C})$,
which is the completion of the linear space $L^0(G_{\alpha },{\bf
C})$ of all simple functions
$$ f(x)=\sum_{j=1}^n b_j \chi _{F_j}(x),$$ where
$b_j\in {\bf C}$, $~F_j\in {\cal A}(G_{\alpha })$, $~F_j\cap
F_k=\emptyset $ for each $j\ne k$, $ ~ \chi _F$ denotes the
characteristic function of a subset $F$, $~\chi _F(x)=1$ for each
$x\in F$ and $\chi _F(x)=0$ for every $x\in G_{\alpha }\setminus F$,
$~n\in {\bf N}$. A norm on $L^1_{G_{\beta }}(G_{\alpha })$ is by our
definition given by the formula:
$$(1)\quad \| f \|_{L^1_{G_{\beta }}(G_{\alpha })} := \sup_{h\in
\theta ^{\beta }_{\alpha }(G_{\beta })} \| f_h \| _{L^1(G_{\alpha
})} <\infty ,$$ where $f_h(g):=f(h^{-1}g)$ for $h, g \in G_{\alpha
}$, $L^1(G_{\alpha }, \mu _{\alpha },{\bf C})$ is the usual Banach
space of all $\mu _{\alpha }$-measurable functions $u: G_{\alpha
}\to {\bf C}$ such that
$$(2)\quad \| u \| _{L^1(G_{\alpha
})} = \int_{G_{\alpha }} |u(g)|\mu _{\alpha }(dg) <\infty .$$
Suppose that \par $(3)$ $\phi : \Lambda \to \Lambda $ is an
increasing mapping, $\alpha < \phi (\alpha ) $ for each $\alpha \in
\Lambda $. We consider the complex space \par $(4)$ $L^{\infty
}(L^1_{G_{\beta }}(G_{\alpha },\mu_{\alpha }): \alpha <\beta \in
\Lambda ) := \{ f=(f_{\alpha }: \alpha \in \Lambda ); ~ f_{\alpha
}\in L^1_{G_{\beta }}(G_{\alpha },\mu_{\alpha })$ $\mbox{ for each
}$ $\alpha \in \Lambda ;$ $\| f \| _{\infty } := \sup_{\alpha \in
\Lambda } \| f_{\alpha } \|_{L^1_{G_{\beta }}(G_{\alpha })} <\infty
,$ $\mbox{ where }$ $\beta = \phi (\alpha ) \} .$ \par When measures
$\mu _{\alpha }$ are specified, spaces are denoted shortly by
$L^1_{G_{\beta }}(G_{\alpha })$ and $L^{\infty }(L^1_{G_{\beta
}}(G_{\alpha }): \alpha <\beta \in \Lambda )$.
\par {\bf 4. Proposition.} {\it Supply the family $L^{\infty
}(L^1_{G_{\beta }}(G_{\alpha }): \alpha <\beta \in \Lambda )$ from
\S 3 with the multiplication $f{\tilde \star } u =w$ such that
$$(1)\quad w_{\alpha } = f_{\beta } {\tilde *} u_{\alpha } =
\int_{G_{\beta }} f_{\beta }(h) u_{\alpha }(\theta ^{\beta }_{\alpha
}(h)g) \mu_{\beta }(dh) .$$
 Then $L^{\infty
}(L^1_{G_{\beta }}(G_{\alpha }): \alpha <\beta \in \Lambda )$ is the
complex normed algebra generally noncommutative and nonassociative.}
\par {\bf Proof.}  Evidently $L^{\infty
}(L^1_{G_{\beta }}(G_{\alpha }): \alpha <\beta \in \Lambda )$ is the
complex linear normed space, since $$\| af+bu \| _{\infty } :=
\sup_{\alpha \in \Lambda } \| af_{\alpha } + bu_{\alpha }
\|_{L^1_{G_{\beta }}(G_{\alpha })} \le \sup_{\alpha \in \Lambda }
|a| \| f_{\alpha } \| _{L^1_{G_{\beta }}(G_{\alpha })} +
\sup_{\alpha \in \Lambda } |b| \| u_{\alpha } \| _{L^1_{G_{\beta
}}(G_{\alpha })}$$ $$ = |a| \| f \| _{\infty } + |b| \| u \|
_{\infty }
$$ for each functions $f, u \in L^{\infty }(L^1_{G_{\beta }}(G_{\alpha }):
\alpha <\beta \in \Lambda )$ and complex numbers $a, b\in {\bf C}$,
where $\beta = \phi (\alpha )$ for each $\alpha \in \Lambda $. On
the other hand,
$$ \| f_{\beta } {\tilde *} u_{\alpha } \|_{L^1_{G_{\beta }}(G_{\alpha })}
\le \| f_{\beta } \|_{L^1(G_{\beta })} \| u_{\alpha }
\|_{L^1_{G_{\beta }}(G_{\alpha })} \le \| f_{\beta }
\|_{L^1_{G_{\beta }}(G_{\alpha })} \| u_{\alpha } \|_{L^1_{G_{\beta
}}(G_{\alpha })}$$ in accordance with Lemma 17.2
\cite{lujms150:4:08}, since $\| f_{\beta } \|_{L^1(G_{\beta })}\le
\| f_{\beta } \|_{L^1_{G_{\beta }}(G_{\alpha })}$, consequently,
$$ \| f{\tilde \star } u \| _{\infty } \le \| f \| _{\infty } \| u \|
_{\infty }.$$
\par The noncommutativity and the nonassociativity of this
multiplication follows from Formula $(1)$.
\par {\bf 5. Corollary.} {\it The family of all nonnegative functions
${\cal P} := \{ f: f \in L^{\infty }(L^1_{G_{\beta }}(G_{\alpha }):
\alpha <\beta \in \Lambda ); ~ f_{\alpha }(x)\ge 0 ~ \forall x\in
G_{\alpha } ~ \forall \alpha \in \Lambda \} $ is a lattice.}
\par {\bf Proof.} If $f, u \in {\cal P}$ and $a\ge 0$ and $b\ge 0$,
then $af+bu\in {\cal P}$, $f\wedge u =\min (f,u)\in {\cal P}$,
$f\vee g=\max (f,g)\in {\cal P}$, since $af_{\alpha }(x)+bu_{\alpha
}(x)\ge 0$ and $\min (f_{\alpha }(x),u_{\alpha }(x))\ge 0$ and $\max
(f_{\alpha }(x),u_{\alpha }(x))\ge 0$ for each $x\in G_{\alpha }$
and $\alpha \in \Lambda $. Moreover, from $\mu _{\alpha }(S)\ge 0$
for each $S\in {\cal B}(G_{\alpha })$ and $\alpha \in \Lambda $ and
from Formula 4$(1)$ it follows that $(f{\tilde \star }u)_{\alpha
}(x)\ge 0$ for each $x\in G_{\alpha }$ and $\alpha \in \Lambda $,
consequently, $f{\tilde \star }u\in {\cal P}$ for all $f, u \in
{\cal P}$.
\par {\bf 6. Corollary.} {\it The operators
$T_f$ and ${\tilde T}_f$ defined by the formulas \par $(1)$ $T_fu:=
f{\tilde \star }u$ and
\par $(2)$ ${\tilde T}_fu:=u{\tilde \star }f$ are $\bf C$-linear and
continuous on the algebra $L^{\infty }(L^1_{G_{\beta }}(G_{\alpha
}): \alpha <\beta \in \Lambda )$ for each $f\in L^{\infty
}(L^1_{G_{\beta }}(G_{\alpha }): \alpha <\beta \in \Lambda )$.}
\par {\bf 7. Lemma.} {\it Let an operator ${\hat U}_g$ be given by the formula:
\par $(1)$ ${\hat U}_{g_{\beta }} f_{\alpha }(x) = f_{\alpha
}(\theta ^{\beta }_{\alpha }(g_{\beta })x)$ for each $\alpha \in
\Lambda $ with $\beta = \phi (\alpha )$, $ ~ g_{\beta }\in G_{\beta
}$, $ ~ x\in G_{\alpha }$, then ${\hat U}_g: L^{\infty
}(L^1_{G_{\beta }}(G_{\alpha }): \alpha <\beta \in \Lambda )\to
L^{\infty }(L^1_{G_{\beta }}(G_{\alpha }): \alpha <\beta \in \Lambda
)$ is the linear isometry for each $g= \{ g_{\beta }: ~ \beta \in
\Lambda , ~ g_{\beta }\in G_{\beta } \} $.}
\par {\bf Proof.} Evidently, ${\hat U}_g (af+bu) = a  {\hat U}_g f +
b {\hat U}_g u$ for each $a, b\in {\bf C}$ and $f, u \in L^{\infty
}(L^1_{G_{\beta }}(G_{\alpha }): \alpha <\beta \in \Lambda )$, since
 ${\hat U}_{g_{\beta }} (af_{\alpha }(x)+bu_{\alpha }(x)) = af_{\alpha
}(\theta ^{\beta }_{\alpha }(g_{\beta })x)+ b f_{\alpha }(\theta
^{\beta }_{\alpha }(g_{\beta })x)= a {\hat U}_{g_{\beta }} f_{\alpha
}(x) + b {\hat U}_{g_{\beta }} u_{\alpha }(x)$ for each $\alpha \in
\Lambda $ with $\beta = \phi (\alpha )$, $ ~ g_{\beta }\in G_{\beta
}$, $ ~ x\in G_{\alpha }$. The isometry property follows from the
equalities:
$$\| {\hat U}_gf \| _{\infty } = \sup _{\alpha \in \Lambda , \beta =\phi (\alpha )}
\sup_{h\in \theta ^{\beta }_{\alpha }(G_{\beta })} \int_{G_{\alpha
}} |f_{\alpha }(h\theta ^{\beta }_{\alpha }(g_{\beta })x)| \mu
_{\alpha }(dx)$$ $$ = \sup _{\alpha \in \Lambda , \beta =\phi
(\alpha )} \sup_{h\in \theta ^{\beta }_{\alpha }(G_{\beta })}
\int_{G_{\alpha }} |f_{\alpha }(hx)| \mu _{\alpha }(dx) = \| f \|
_{\infty } .$$
\par {\bf 8. Lemma.} {\it The operator ${\hat U}_g$ from \S 7
satisfies the equality:
$$(1)\quad {\hat U}_g(f{\tilde \star }u) = f {\tilde \star }{\hat
U}_gu$$ for each $f, u \in L^{\infty }(L^1_{G_{\beta }}(G_{\alpha
}): \alpha <\beta \in \Lambda )$ and every $g= \{ g_{\beta }: ~
\beta \in \Lambda , ~ g_{\beta }\in G_{\beta } \} $.}
\par {\bf Proof.} Formula $(1)$ follows from the equalities:
$${\hat U}_{g_{\beta }}(f_{\beta } {\tilde *}u_{\alpha })(x) =
\int_{G_{\beta }} f_{\beta }(h)u_{\alpha }((\theta ^{\beta }_{\alpha
}(hg_{\beta })x))\mu _{\beta }(dh) = (f _{\beta } {\tilde
*} ({\hat U}_{g_{\beta }}u_{\alpha }))(x)$$ for each $x\in G_{\alpha
}$, $~\alpha \in \Lambda $, $~\beta = \phi (\alpha )$.
\par {\bf 9. Corollary.} {\it A bijective correspondence
between elements $g\in \prod_{\alpha \in \Lambda }G_{\alpha }=: G$
and operators ${\hat U}_g$ exists so that \par $(1)$ ${\hat
U}_g{\hat U}_h={\hat U}_{gh}$ for each $g, h\in G$.}
\par {\bf Proof.} The group $G$ consists of elements
$g= \{ g_{\alpha }: ~ \alpha \in \Lambda , ~ g_{\alpha }\in
G_{\alpha } \} $ with the multiplication $gh = \{ g_{\alpha
}h_{\alpha }: ~ \alpha \in \Lambda \} $ and the inversion $g^{-1} =
\{ g_{\alpha }^{-1}: \alpha \in \Lambda \} $, since $g_{\alpha
}h_{\alpha }\in G_{\alpha } $ for each $\alpha \in \Lambda $. This
group $G$ is the topological group relative to the Tychonoff
(product) topology $\tau ^t$ with the base $V=V_{\alpha }\times
\prod_{\gamma \in \Lambda ; \gamma \ne \alpha } G_{\gamma }$, where
$V_{\alpha }\in \tau _{\alpha }$, $ ~ \tau _{\alpha }$ is the
topology on $G_{\alpha }$, $~ \alpha \in \Lambda $. It is also the
topological group relative to the box topology $\tau ^b$ with the
base $V=\prod_{\alpha \in \Lambda } V_{\alpha }$, where $V_{\alpha
}\in \tau _{\alpha }$ is open in $G_{\alpha }$ for each $\alpha \in
\Lambda $. Then we deduce that
$${\hat U}_{g_{\beta }}{\hat U}_{h_{\beta }}f_{\alpha }(x)={\hat
U}_{g_{\beta }} ({\hat U}_{h_{\beta }} f_{\alpha }(x)) = {\hat
U}_{g_{\beta }}(f_{\alpha }(\theta ^{\beta }_{\alpha }(h_{\beta
})x))=$$ $$f_{\alpha } (\theta ^{\beta }_{\alpha }(g_{\beta
})(\theta ^{\beta }_{\alpha }(h_{\beta })x))=f_{\alpha } (\theta
^{\beta }_{\alpha }(g_{\beta }h_{\beta })x)= {\hat U}_{g_{\beta
}h_{\beta }}f_{\alpha }(x)$$ for each $\alpha \in \Lambda $ with
$\beta = \phi (\alpha )$, where $x\in G_{\alpha }$. The latter
relation implies Formula $(1)$.
\par {\bf 10. Proposition.} {\it The representation ${\hat U}:
(G,\tau ^b) \to L^{\infty }(L^1_{G_{\beta }}(G_{\alpha }): \alpha
<\beta \in \Lambda )$ is strongly continuous.}
\par {\bf Proof.} Each bounded either continuous or simple function $f: G_{\alpha }\to
{\bf C}$ evidently belongs to $L^1_{G_{\beta }}(G_{\alpha })$, since
$\mu _{\alpha }$ is the probability measure and
$$(1)\quad \| f \|_{L^1_{G_{\beta }}(G_{\alpha })} \le \sup_{x\in G_{\alpha
}} |f(x)|<\infty .$$ Each compact subset $V$ in $G_{\alpha }$ is
closed in $G_{\alpha }$ in accordance with Theorem 3.1.8 \cite{eng},
consequently, every compact subset $V$ is a Borel subset. \par  For
an arbitrary marked function $u\in L^1_{G_{\beta }}(G_{\alpha })$
from the inclusion $L^1_{G_{\beta }}(G_{\alpha })\subset
L^1(G_{\alpha })$, the Borel regularity of the measure $\mu _{\alpha
}$, Conditions 2$(1-4)$ and Lusin's theorem 2.3.5 \cite{federer} it
follows that for each $\epsilon _n>0$ a compact subset $E_n=E_n(u)$
in $G_{\alpha }$ exists so that the restriction $u|_{E_n}$ is
continuous and $\mu _{\alpha }(G_{\alpha }\setminus E_n)<\epsilon
_n$, since the measure $\mu _{\alpha }$ is radonian and for each
$\delta >0$ a compact subset $V$ in $G_{\alpha }$ exists such that
$\mu _{\alpha }(G_{\alpha }\setminus V)<\delta $ and considering
$u|_V$ and $\mu _{\alpha }|_{{\cal B}(V)}$. Take a monotone
decreasing sequence $\epsilon _n$ such that $\epsilon _n\downarrow
0$. Then $\mu _{\alpha }(G_{\alpha }\setminus E)=0$, where $$E(u)= E
:= \bigcup_{n=1}^{\infty }E_n.$$ For every $\delta >0$ and each
restriction $u|_{E_n}$ a simple function $$v_n = \sum_{j=1}^{m(n)}
b_{j,n} \chi _{F_{j,n}}(x)$$ exists such that $$\sup_{x\in E_n}
|u(x)-v_n(x)|<\delta ,$$ where $b_{j,n}\in {\bf C},$ $~ F_{j,n}\in
{\cal B}(E_n)$, $~ m(n)\in {\bf N}$. Put $v_0(x)=0$ on $G_{\alpha
}\setminus E$ and take the combination $v$ of these mappings $v_n$,
then $\sup_{x\in E} |u(x)-v(x)|<\delta $ and hence $v\in
L^1_{G_{\beta }}(G_{\alpha })$ due to Inequality $(1)$. \par In
accordance with \S 3 in each space $L^1_{G_{\beta }}(G_{\alpha })$
the linear space of all simple functions
$$(2)\quad f(x)=\sum_{j=1}^p b_j \chi _{F_j}(x)$$ is dense, where $b_j\in {\bf
C}$, $~F_j\in {\cal A}(G_{\alpha })$, $\chi _F$ denotes the
characteristic function of a subset $F$, $~\chi _F(x)=1$ for each
$x\in F$ and $\chi _F(x)=0$ for every $x\in G_{\alpha }\setminus F$,
$~p\in {\bf N}$. \par In view of Lemma 7 it is sufficient to prove,
that the representation $G\ni g \mapsto {\hat U}_gf_{\alpha }$ is
continuous on each simple function $f_{\alpha }\in L^1_{G_{\beta
}}(G_{\alpha })$, when $G$ is supplied with the box topology $\tau
^b$,  since $\| {\hat U}_g(f_{\alpha }-u) \| = \| f_{\alpha }-u \|
$.
\par Now we take $E_n=E_n(f_{\alpha })$ and $E=E(f_{\alpha })$ as above. The measure
$\mu _{\alpha }$ is quasi-invariant, consequently,
\par $(3)$ $\mu _{\alpha }(h^{-1}(G_{\alpha }\setminus E))=0$ and
hence $\mu _{\alpha }([h^{-1}E]\cap W_{\alpha })=1$ for each $h\in
\theta ^{\beta }_{\alpha }(G_{\beta })$ with $\beta = \phi (\alpha
)$ and $\alpha \in \Lambda $.
\par On the other hand,
\par $(4)$ ${\hat U}_{g_{\beta }} f_{\alpha }(x) - {\hat U}_{e_{\beta }}f_{\alpha }(x) =
f_{\alpha }(\theta ^{\beta }_{\alpha }(g_{\beta })x) - f_{\alpha
}(x)$ and
$$ \| [{\hat U}_{g_{\beta }}  - {\hat U}_{e_{\beta }}] f_{\alpha }(x)\|_{ L^1_{G_{\beta }}(G_{\alpha })} =
\sup_{h\in \theta ^{\beta }_{\alpha }(G_{\beta })} \int_{G_{\alpha
}} |f(h\theta ^{\beta }_{\alpha }(g_{\beta })x) - f(hx)| \mu
_{\alpha }(dx) \le 2 \| f \|_{L^1_{G_{\beta }}(G_{\alpha })}<\infty
$$ for each $\alpha \in \Lambda $ with $\beta = \phi (\alpha )$, $ ~
g_{\beta }\in G_{\beta }$, $ ~ x\in G_{\alpha }$, since $\theta
^{\beta }_{\alpha }(e_{\beta })=e_{\alpha }$ is the unit element in
the group $G_{\alpha }$. \par For each $\delta >0$ an element
$h_{\delta }\in \theta ^{\beta }_{\alpha }(G_{\beta })$ exists such
that
$$(5)\quad  | \| f_{\alpha }(x)\|_{L^1_{G_{\beta }}(G_{\alpha })} -
\int_{G_{\alpha }} |f_{\alpha }(h_{\delta }x)| \mu _{\alpha }(dx) |
<\delta .$$ Evidently, the series $$(6)\quad \int_{G_{\alpha }}
|f(hx)| \mu _{\alpha }(dx)=\sum_{j=1}^p |b_j| \mu _{\alpha
}(h^{-1}F_j)\le \sum_{j=1}^p |b_j|$$ is finite, since $\mu _{\alpha
}$ is the nonnegative measure and $\mu _{\alpha }(G_{\alpha })=1$
and $p\in {\bf N}$ is a natural number. Each measure
$$(7)\quad w_h(A):= \sum_{j=1}^p |b_j| \mu _{\alpha }
((h^{-1}F_j)\cap A)$$ is $\sigma $-additive on ${\cal A}(G_{\alpha
})$ and absolutely continuous relative to $\mu _{\alpha }$ due
Formula $(6)$, where $A\in {\cal A}(G_{\alpha })$, $~h\in G_{\beta
}$. \par For a given arbitrary positive number $\epsilon >0$ take a
natural number $n_0$ such that $\epsilon >\epsilon _n$ for each
$n>n_0$. Choose a marked natural number $m>n_0$. From Theorem 4.5
\cite{hew} and Formulas $(1-4,7)$ and Conditions 1$(1-4)$ and
2$(1-4)$ it follows that for each $\delta >0$ a symmetric $V_{\beta
}=V_{\beta }^{-1}$ neighborhood of $e_{\beta }$ in $G_{\beta }$
exists such that $\theta ^{\beta }_{\alpha }(V_{\beta })$ is
precompact in $G_{\alpha }$ and
$$(8)\quad \int_{G_{\alpha }\cap E_m}
|\sum_{j=1}^p b_j \chi _{F_j}(\theta ^{\beta }_{\alpha }(g_{\beta
})hx)- \sum_{j=1}^p b_j \chi _{F_j}(hx) |\mu _{\alpha }(dx)<\delta
.$$ This is possible by a choice of a sufficiently small $V_{\beta
}$ such that the left quasi-invariance factor $\rho ^l(a,b)$ is
bounded on $V_{\beta }\times (W_{\alpha }\cap E_m)$, since $\rho _l$
is continuous on $G_{\beta }\times W_{\alpha }$ and
$E_m=E_m(f_{\alpha })$ is compact. Indeed, the product $cl(\theta
^{\beta }_{\alpha }(V_{\beta }))E_m=:Q_m$ is compact in $G_{\alpha
}$ as the product of two compact subsets in the topological group
$G_{\alpha }$ (see \S 4.4 in \cite{hew}) and $E_m\subset Q_m$. From
the choice of $E_m$ we infer, that
$$(9)\quad \int_{G_{\alpha }\setminus E_m} |\sum_{j=1}^p b_j \chi
_{F_j}(\theta ^{\beta }_{\alpha }(g_{\beta })hx)- \sum_{j=1}^p b_j
\chi _{F_j}(hx) |\mu _{\alpha }(dx)\le 2 \epsilon \sum_{j=1}^p |b_j|
$$ for each $h\in \theta ^{\beta }_{\alpha }(G_{\beta })$. \par Thus
from $(8,9)$ it follows ,that for each $\delta >0$ a neighborhood
$V_{\beta }$ of $e_{\beta }$ in $G_{\beta }$ exists such that $\|
[{\hat U}_{g_{\beta }}  - {\hat U}_{e_{\beta }}] f_{\alpha }(x)\|_{
L^1_{G_{\beta }}(G_{\alpha })}<\delta $ for each $g_{\beta }\in
V_{\beta }$. Taking $V=\prod_{\alpha \in \Lambda } V_{\alpha }$ we
get that $\| [{\hat U}_g  - {\hat U}_e] f \| <\delta $ for every
$g\in V$, where $f=(f_{\alpha }: \alpha \in \Lambda )\in L^{\infty
}(L^1_{G_{\beta }}(G_{\alpha }): \alpha <\beta \in \Lambda )$.
\par {\bf 11. Theorem.} {\it Suppose that a continuous mapping $S: L^{\infty
}(L^1_{G_{\beta }}(G_{\alpha }): \alpha <\beta \in \Lambda )\to
L^{\infty }(L^1_{G_{\beta }}(G_{\alpha }): \alpha <\beta \in \Lambda
)$ satisfies the following conditions:
\par $(1)$ $S$ is linear over the complex field so that $Sf=(S_{\alpha }f_{\alpha }: \alpha \in \Lambda )$
with $S_{\alpha }f_{\alpha }\in L^1_{G_{\beta }}(G_{\alpha })$ for
each $\alpha \in \Lambda $, where $\beta =\phi (\alpha )$;
\par $(2)$ positive, i.e. $S_{\alpha }f_{\alpha }$ is positive
if $f_{\alpha }$ is positive for each $\alpha \in \Lambda $;
\par $(3)$ $S(f{\tilde \star }g)=f{\tilde \star }Sg$ for every $f,
g\in L^{\infty }(L^1_{G_{\beta }}(G_{\alpha }): \alpha <\beta \in
\Lambda )$. \par Then elements $a\in G$ and $p= \{ p_{\alpha }:
p_{\alpha }>0 ~ \forall \alpha \in \Lambda \} \in {\bf R}^{\Lambda
}$ exist so that
\par $(4)$ $S=p{\hat U}_a$, that is $S_{\alpha }f_{\alpha }(x) =
p_{\alpha }{\hat U}_{a_{\beta }}f_{\alpha }(x)$ for any $\alpha \in
\Lambda $ with $\beta =\phi (\alpha )$ and each $x\in G_{\alpha }$.}
\par {\bf Proof.} In view of Corollary 5 the family ${\cal P}$ of all
nonnegative elements forms the cone. Conditions $(1,2)$ imply that
\par $(5)$ $S(f\vee g)=(Sf)\vee (Sg)$ and $S(f\wedge g)=(Sf)\wedge
(Sg)$ and
\par $(6)$ $S(\vee_{n=1}^{\infty } g_n)=\vee_{n=1}^{\infty } (Sg_n)$ and
$S(\wedge_{n=1}^{\infty } g_n)=\wedge_{n=1}^{\infty } (Sg_n)$ on
${\cal P}$. \par Being continuous the operator $S$ is bounded. We
consider a subset $E=\prod_{\alpha \in \Lambda } E_{\alpha }$ such
that $E_{\alpha }\in {\cal A}(G_{\alpha })$ for each $\alpha $. The
function \par $(7)$ $\xi _E = \{ \chi _{E_{\alpha }}: \alpha \in
\Lambda \} $ belongs to $L^{\infty }(L^1_{G_{\beta }}(G_{\alpha }):
\alpha <\beta \in \Lambda )$, since $\int_{G_{\alpha }} \chi
_{E_{\alpha }}(x)\mu _{\alpha }(dx) = \mu _{\alpha }(E_{\alpha })\le
\mu _{\alpha }(G_{\alpha })=1$ for any $\alpha $. The function $\eta
_E (x) := S\xi _E(x)$ is positive by the conditions of this theorem,
hence $\eta _{E,\alpha } := S_{\alpha } \chi _{E_{\alpha
}}(x_{\alpha })\ge 0$ on $G_{\alpha }$ and \par $(8)$ $\eta
_{E,\alpha }(x_{\alpha })>0$ for each $x_{\alpha }\in T_{S,\alpha
}(E_{\alpha })$ and every $\alpha \in \Lambda $, where $T_{S,\alpha
}(E_{\alpha })$ denotes a subset in $G_{\alpha }$  on which a
function $\eta _{E,\alpha }$ is positive which is defined up to a
$\mu _{\alpha }$-null set. From Formulas $(5)$ we get that
\par $(9)$ $T_{S,\alpha }(E_{\alpha }\cup F_{\alpha })=T_{S,\alpha }(E_{\alpha
})\cup T_{S,\alpha }(F_{\alpha })$ and $T_{S,\alpha }(E_{\alpha
}\cap F_{\alpha })=T_{S,\alpha }(E_{\alpha })\cap T_{S,\alpha
}(F_{\alpha })$ for each $\alpha \in \Lambda $ and $E_{\alpha },
F_{\alpha }\in {\cal A}(G_{\alpha })$. Moreover, the definition of
$T_{S,\alpha }(E_{\alpha })$ by Formula $(8)$ implies that \par
$(10)$ $(T_S)^{-1}=T_{S^{-1}}$, that is $(T_{S,\alpha
})^{-1}=T_{S^{-1},\alpha }$ for any $\alpha \in \Lambda $.
\par For a marked $\alpha \in \Lambda $ and $\beta = \phi (\alpha )$ we next
consider a base of symmetric neighborhoods $U_{\alpha
,v}=U^{-1}_{\alpha ,v}$ and $U_{\beta ,v}=U^{-1}_{\beta ,v}$ of the
unit elements $e_{\alpha }$ in $G_{\alpha }$ and $e_{\beta }$ in
$G_{\beta }$ satisfying the conditions:
\par $(11)$ $U_{\beta ,v}\subseteq U_{\alpha ,v}\cap G_{\beta }$ for
each $v\in \Upsilon _{\alpha ,\beta }$, where \par $(12)$ $\Upsilon
_{\alpha ,\beta }$ is a directed set by inclusion: $v\le t$ if and
only if $U_{\alpha ,t}\subseteq U_{\alpha ,v}$ so that for each
$v\in \Upsilon _{\alpha ,\beta }$ there exists $q(v)\in \Upsilon
_{\alpha ,\beta }$ with $q(v)>v$ and $U_{\alpha ,q(v)} \subset
U_{\alpha ,v}$ and $U_{\alpha ,q(v)} \ne U_{\alpha ,v}$,
$~U^{-1}_{\alpha ,v} = \{ g^{-1}: ~ g \in U_{\alpha ,v} \} $. Then
we put
\par $(13)$ $\xi _{\alpha ,v} = \chi _{U_{\alpha ,v}}$, $\eta
_{\alpha ,v}=S_{\alpha }\xi _{\alpha ,v}$, $\xi _v := \{ \xi
_{\alpha ,v}: \alpha \in \Lambda \} $ and $\eta _v := \{ \eta
_{\alpha ,v}: \alpha \in \Lambda \} $ for any $v\in \Upsilon
_{\alpha ,\beta }$. Below the proof of this theorem is continued and
is based on the following intermediate lemmas.
\par {\bf 12. Lemma.} {\it Let $P=T_S^{-1}U$, that is $P_{\alpha ,v} = T_{S,\alpha }^{-1}(U_{\alpha ,v})$
for each $\alpha \in \Lambda $ and $v\in \Upsilon _{\alpha ,\beta
}$, where $T_S$ and $U_{\alpha ,v}$ and $\Upsilon _{\alpha ,\beta }$
with $\beta =\phi (\alpha )$ are as in \S 11. Then for any $\alpha
\in \Lambda $ and for each $v\in \Upsilon _{\alpha ,\beta }$
elements $w=w(v)\in \Upsilon _{\alpha ,\beta }$ and $a_{\beta ,v}\in
G_{\beta }$ exist such that
\par $(1)$ $a_{\beta ,v} U_{\beta ,v}\supset P_{\beta ,w(v)}$
up to a $\mu _{\beta }$-null set.}
\par {\bf Proof.} Suppose the contrary that there exists $U_{\beta
,v}$ so that $a_{\beta }U_{\beta ,v}$ does not cover $\mu _{\beta
}$-almost entirely the set $P_{\beta ,w}$ for any $w\in \Upsilon
_{\alpha ,\beta }$ and any $a_{\beta }\in G_{\beta }$. Since
$U_{\alpha ,v}$ is a base of neighborhoods of the unit element in
$G_{\alpha }$ there exist $t, s \in \Upsilon _{\alpha ,\beta }$ such
that $U_{\beta ,t}^3\subseteq U_{\beta ,v}$ and $U_{\beta
,s}^2\subseteq U_{\beta ,v}$. Take two elements $a_{\beta ,1},
a_{\beta ,2}\in G_{\beta }$ satisfying the conditions:
\par $(2)$ $A_{\alpha } =cl_{\alpha } \theta ^{\beta }_{\alpha }([a_{\beta ,1}U_{\beta ,t}]
\cap P_{\beta ,s})$ and $B_{\alpha } =cl_{\alpha } \theta ^{\beta
}_{\alpha }([a_{\beta ,2}U_{\beta ,t}]\cap P_{\beta ,s})$ and
$A_{\beta } =[a_{\beta ,1} U_{\beta ,t}]\cap P_{\beta ,s}$ and
$B_{\beta } =[a_{\beta ,2}U_{\beta ,t}]\cap P_{\beta ,s}$ with $\mu
_{\alpha }(A_{\alpha })>0$ and $\mu _{\alpha }(B_{\alpha })>0$ and
$\mu _{\beta }(A_{\beta })>0$ and $\mu _{\beta }(B_{\beta })>0$ and
\par $(3)$ $A_{\beta }\cap B_{\beta }U_{\beta ,t}=\emptyset $, where
$cl_{\alpha }(A)$ denotes the closure of a subset $A$ in $G_{\alpha
}$. This is possible, since the group $\theta ^{\beta }_{\alpha
}(G_{\beta })$ is dense in $G_{\alpha }$ and the quasi-invariant
radonian Borel regular  measures $\mu _{\alpha }$ and $\mu _{\beta
}$ are positive on each open subset in $G_{\alpha }$ and $G_{\beta
}$ correspondingly.
\par Then the sets \par $(4)$ $C_{\alpha }=T_{S,\alpha }^{-1}(A_{\alpha })$
and $D_{\alpha }=T_{S,\alpha }^{-1}(B_{\alpha })$ are $\mu _{\alpha
}$-measurable and $~C_{\alpha }\cup D_{\alpha } \subset cl_{\alpha
}(U_{\beta ,s})$, also $C_{\beta }=T_{S,\beta }^{-1}(A_{\beta })$
and $D_{\beta }=T_{S,\beta }^{-1}(B_{\beta })$ are $\mu _{\beta
}$-measurable and are contained in $U_{\beta ,s}$, $~C_{\beta }\cup
D_{\beta } \subset U_{\beta ,s}$. From Conditions 11$(10,11)$ we
deduce that \par $(5)$ $C_{\alpha }^{-1}D_{\alpha }\subset
cl_{\alpha }(U_{\beta , v})$ and $C_{\beta }^{-1}D_{\beta }\subset
U_{\beta , v}$. \par Applying Formula 11$(3)$ we infer the
following:
$$(6)\quad (\chi_{C_{\beta }} {\tilde *} \chi_{cl_{\alpha }(U_{\beta , t})})(x)
= \int_{G_{\beta }} \chi_{C_{\beta }}(h) \chi_{cl_{\alpha }(U_{\beta
, t})}(\theta ^{\beta }_{\alpha }(h)x) \mu _{\beta }(dh)$$  $$ =
\int_{C_{\beta }} \chi_{cl_{\alpha }(U_{\beta , t})}(\theta ^{\beta
}_{\alpha }(h)x) \mu _{\beta }(dh) = \mu _{\beta } (C_{\beta }\cap
cl_{\alpha }(U_{\beta ,t}x^{-1})) ,$$ consequently, $(\chi_{C_{\beta
}} {\tilde *} \chi_{cl_{\alpha }(U_{\beta , t})})(x)=\mu_{\beta
}(C_{\beta })>0$ for each $x\in D_{\alpha }$, since $C_{\alpha
}^{-1}D_{\alpha } \subset cl_{\alpha }(U_{\beta ,t})$ and $C_{\alpha
}D_{\alpha }^{-1} \subset cl_{\alpha }(U_{\beta ,t}^{-1})=cl_{\alpha
}(U_{\beta ,t})$ and hence
$$(7)\quad (\chi_{C_{\beta }} {\tilde *} \chi_{cl_{\alpha }(U_{\alpha ,
t})})(x)\ge \mu_{\beta }(C_{\beta })\chi_{D_{\alpha }}(x).$$ Then
the inequality
$$(8) \quad S(\chi_{C_{\beta }} {\tilde *} \chi_{cl_{\alpha }(U_{\beta ,
t})})  = (\chi_{C_{\beta }} {\tilde *} S\chi_{cl_{\alpha }(U_{\beta
, t})}) \ge \mu_{\beta }(C_{\beta })S\chi_{D_{\alpha }}$$ follows
from Formulas $(7)$ and 11$(3)$. Applying Conditions $(2)$ one gets
$A_{\alpha }=T_{S,\alpha } (C_{\alpha }) = \{ x: ~ S\chi _{C_{\alpha
}} (x)>0 \} $ and $B_{\alpha }=T_{S,\alpha } (D_{\alpha }) = \{ x: ~
S\chi _{D_{\alpha }} (x)>0 \} $ and $A_{\beta }=T_{S,\beta }
(C_{\beta }) = \{ x: ~ S\chi _{C_{\beta }} (x)>0 \} $ and $B_{\beta
}=T_{S,\beta } (D_{\beta }) = \{ x: ~ S\chi _{D_{\beta }} (x)>0 \}
$. On the other hand, from Formulas $(3)$ we deduce the following:
$$(9)\quad S(\chi_{C_{\beta }} {\tilde *} \chi_{cl_{\alpha }(U_{\beta ,
t})})(x)= S\int_{G_{\beta }\cap U_{\beta, t}x^{-1}} \chi_{C_{\beta
}}(h) \mu _{\beta }(dh)=0,$$ since $A_{\beta }\cap (B_{\beta
}U_{\beta ,t})=\emptyset $ and hence $T_{S,\beta }^{-1}(A_{\beta
})\cap T_{S,\beta }^{-1}(B_{\beta }U_{\beta ,t})=\emptyset $ so that
$C_{\beta }\cap (D_{\beta }P_{\beta ,t}) = \emptyset $. But Formula
$(9)$ contradicts $(7)$, that finishes the proof of this lemma.
\par {\bf 13. Lemma.} {\it Let an
operation ${\tilde *}$ on $M(G_{\beta })\times L^1_{G_{\beta
}}(G_{\alpha })$ be defined by the formula
$$(1)\quad (\nu {\tilde *} u)(g) = \int_{G_{\beta }} \nu (dh) u(\theta ^{\beta }_{\alpha }(h)g)$$
for each $g \in G_{\alpha }$, where $M(G_{\beta })$ denotes the set
of all finite radon measures $\nu $ on $G_{\beta }$ supplied with
the norm. Then the mapping ${\tilde *}: M(G_{\beta })\times
L^1_{G_{\beta }}(G_{\alpha })\to L^1_{G_{\beta }}(G_{\alpha })$ is
bilinear and continuous.}
\par {\bf Proof.} From Formula $(1)$ it follows that
$(\nu {\tilde *} (au+bg)) = a(\nu {\tilde *} u)+b(\nu {\tilde *} g)$
and $(a\nu +b\mu ) {\tilde *} u = a(\nu {\tilde *} u) +b (\mu
{\tilde *} u)$ for any complex numbers $a, b \in {\bf C}$, radonian
measures $\nu , \mu \in M(G_{\beta })$ and functions $u, g \in
L^1_{G_{\beta }}(G_{\alpha })$. Remind that the space $M(G_{\beta
})$ is supplied with the standard norm:  $\| \nu \| = |\nu
|(G_{\beta })$, where $|\nu |=\nu ^+ + \nu ^-$ is the variation of
$\nu $ with the standard decomposition $\nu =\nu ^+ - \nu ^-$ into
the difference of two nonnegative measures $\nu ^+$ and $\nu ^-$.
Then $$\sup_{s\in \theta ^{\beta }_{\alpha }(G_{\beta
})}\int_{G_{\alpha }} |\int_{G_{\beta }} \nu (dh) u(\theta ^{\beta
}_{\alpha }(h)(sg))|\mu_{\alpha }(dg)\le $$  $$\sup_{s\in \theta
^{\beta }_{\alpha }(G_{\beta })}\int_{G_{\beta }} \{ |\nu |(dh)
\int_{G_{\alpha }} |u(\theta ^{\beta }_{\alpha }(h)(sg))|\mu_{\alpha
}(dg) \} $$ due to Fubini's theorem, consequently,
$$(2)\quad \| \nu {\tilde *} u \|_{L^1_{G_{\beta }}(G_{\alpha
})}\le \| \nu \| \| u \|_{L^1_{G_{\beta }}(G_{\alpha })}.$$
 The latter inequality implies the continuity
of such skew convolution ${\tilde *}: M(G_{\beta })\times
L^1_{G_{\beta }}(G_{\alpha })\to L^1_{G_{\beta }}(G_{\alpha })$.

\par {\bf Continuation of the proof of Theorem 11.}
By transfinite induction and Teichm\"uller-Tukey's lemma applying
Lemma 12 one gets a base of symmetric neighborhoods of the unit
elements such that
\par $(14)$ $U_{\beta ,t}^3\subseteq U_{\beta ,v}$ and $U_{\beta
,s}^2\subseteq U_{\beta ,v}$ for each $t\in \lambda (\beta ,v)$ and
$s\in \nu (\beta ,v)$, where $\lambda (\beta ,v)\subset \nu (\beta
,v)$ are cofinal subsets in $\Upsilon _{\alpha ,\beta }$ all
elements of which are greater than $v$ (see also \S 1.3 \cite{eng});
and Conditions 11$(11,12)$ and 12$(1)$. The inclusion $P_{\beta
,v}\supseteq P_{\beta , q}$ for each $v<q$ leads to $a_{\beta ,v}
U_{\beta ,v} \cap a_{\beta ,q} U_{\beta ,q}\ne \emptyset $,
consequently,
\par $(15)$ $a_{\beta ,v}^{-1} a_{\beta ,q} \in U_{\beta ,v}U_{\beta
,q}^{-1}$ \\ and hence $\{ a_{\beta ,v}: v \} $ is a fundamental
(Cauchy) net in $G_{\beta }$. But $G_{\beta }$ is the topological
group complete relative to its left uniformly  (see \S \S 8.1.17 and
8.3 \cite{eng}). Therefore, this net converges $\lim_v a_{\beta
,v}=a_{\beta }$ in $G_{\beta }$. From $(15)$ the inclusion $a_{\beta
}\in a_{\beta ,v}U_{\beta ,v}^2$ follows, consequently,
\par $(16)$ $a_{\beta } U_{\beta ,v}\supset  a_{\beta ,s} U_{\beta
,s}\supset P_{\beta ,w(s)}$ for every $s\in \nu (\beta ,v)$. \par In
view of Proposition 17.7 \cite{lujms150:4:08}
$$(17)\quad \lim_s \|  \frac{\chi _{U_{\beta ,w(s)}} }{\mu _{\beta }(U_{\beta ,w(s)})}
{\tilde *}  f-f \| _{L^1_{G_{\beta }}(G_{\alpha })}=0$$ for each
$f\in L^1_{G_{\beta }}(G_{\alpha })$, since $\mu _{\beta }(U_{\beta
,w(s)})>0$, where $\beta =\phi (\alpha )$. \par The left
quasi-invariance factor $\rho^l_{\mu _{\alpha }}(\theta ^{\beta
}_{\alpha }(h),g)$ is continuous on $G_{\beta }\times W_{\alpha }$
and satisfies the cocycle condition $$\rho^l_{\mu _{\alpha }}(\theta
^{\beta }_{\alpha }(h),\theta ^{\beta }_{\alpha
}(t^{-1})g)\rho^l_{\mu _{\alpha }}(\theta ^{\beta }_{\alpha
}(t),g)=\rho^l_{\mu _{\alpha }}(\theta ^{\beta }_{\alpha }(ht),g)$$
for all $h, t \in G_{\beta }$ and $g\in G_{\alpha }$. The
probability measure $\mu _{\alpha }$ is Borel regular and $\mu
_{\alpha }(G_{\alpha })=\mu _{\alpha }(W_{\alpha })$, consequently,
$W_{\alpha }$ is dense in $G_{\alpha }$ and hence has a continuous
extension onto $G_{\beta }\times G_{\alpha }$ due to the cocycle
condition and since $\tau _{\alpha }\cap {\theta ^{\beta }_{\alpha
}(G_{\beta })} \subset \theta ^{\beta }_{\alpha }(\tau _{\beta })$,
where $\beta =\phi (\alpha )$. Henceforth, we denote this continuous
extension by the same symbol $\rho^l_{\mu _{\alpha }}$.
\par  Take a net of bounded functions $f_{\alpha ,\kappa
,y}\in L^1_{G_{\beta }}(G_{\alpha })$ such that
$$(18)\quad \lim_{\kappa } \int_{G_{\alpha }}f_{\alpha ,\kappa ,y} (g) \mu
_{\alpha }(dg)\mu _{\beta }(dh)=\delta _y(dh)$$ in $M(G_{\beta })$,
where $\delta _y(dh)$ denotes the atomic Dirac measure on $G_{\beta
}$ with atom at $y\in G_{\beta }$, $\kappa \in K$, where $K$ is a
directed set. Without loss of generality these functions can be
chosen such that the linear span over the complex field $\bf C$ of
the family of functions $span_{\bf C} \{ f_{\alpha ,\kappa ,y}: y
\in G_{\beta }, \kappa _0<\kappa \in K \} $ is dense in
$L^1_{G_{\beta }}(G_{\alpha })$ for each $\kappa _0\in K$. \par A
subspace of continuous functions in $L^1_{G_{\beta }}(G_{\alpha })$
is dense in this space $L^1_{G_{\beta }}(G_{\alpha })$ (see also \S
10). \par In the space $L^1(G_{\beta })$ the net $\frac{\chi
_{U_{\beta ,w(s)}}}{\mu _{\beta }(U_{\beta ,w(s)})}$ converges to
the atomic Dirac measure $\delta _{e_{\beta }}$ on $G_{\beta }$.
Consider these family of functions related by the left shifts with
weight factors $$f_{\alpha ,\kappa ,y} (g)=\frac{f_{\alpha ,\kappa
,e_{\alpha }}(\theta^{\beta }_{\alpha }(y)g)}{\rho^l_{\mu _{\alpha
}} (\theta^{\beta }_{\alpha }(y),e_{\alpha })},$$ where $y\in
G_{\beta }$, then we infer that
$$(19)\quad \int_{G_{\alpha }}f_{\alpha ,\kappa ,e_{\alpha }} (\theta^{\beta }_{\alpha }(hy)g) \mu
_{\alpha }(dg)\mu _{\beta }(dh)= \int_{G_{\alpha }}f_{\alpha ,\kappa
,e} (s) \frac{\rho^l_{\mu _{\alpha }} (\theta^{\beta }_{\alpha
}(hy),s)}{\rho^l_{\mu _{\alpha }} (\theta^{\beta }_{\alpha
}(y),e_{\alpha })} \mu _{\alpha }(ds)\mu _{\beta }(dh)$$ and
$$(20) \quad \lim_s (\frac{\chi _{U_{\beta ,w(s)}}}
{\mu _{\beta }(U_{\beta ,w(s)})} {\tilde *} f_{\alpha ,\kappa
,e_{\alpha }} (s) \frac{\rho^l_{\mu _{\alpha }} (\theta^{\beta
}_{\alpha }(hy),s))}{\rho^l_{\mu _{\alpha }} (\theta^{\beta
}_{\alpha }(y),e_{\alpha })} = f_{\alpha ,\kappa ,e_{\alpha }}
(g)\frac{\rho^l_{\mu _{\alpha }} (\theta^{\beta }_{\alpha
}(y),g)}{\rho^l_{\mu _{\alpha }} (\theta^{\beta }_{\alpha
}(y),e_{\alpha })},$$ since $e_{\beta }g=g$.
\par Consider particularly $f_{\alpha ,v,e_{\alpha }}=\frac{\chi _{U_{\alpha
,v}}} {\mu _{\alpha }(U_{\alpha ,v})}$
  Applying Formulas $(17-20)$ and 8$(1)$ and Lemma 13
we deduce that the limit
$$\lim_s \lim_v (\frac{\chi _{U_{\beta ,w(s)}}}
{\mu _{\beta }(U_{\beta ,w(s)})} {\tilde *} \frac{\xi _{\alpha ,v}}
{\mu _{\alpha }(U_{\alpha ,v})}) = p_{\alpha }>0$$ converges and is
independent of $h\in G_{\beta }$, where $\xi _{\alpha ,v} := \|
S_{\alpha }\chi _{U_{\alpha ,v}} \| $, since $S$ is a bounded linear
operator and
$$ \lim_s \|  ( \frac{\chi _{U_{\beta ,w(s)}} }{\mu _{\beta }(U_{\beta ,w(s)})}
{\tilde *}  S_{\alpha }f) - S_{\alpha }f \| _{L^1_{G_{\beta
}}(G_{\alpha })}=0$$ and
$$\lim_{\kappa }  r {\tilde *}  f_{\alpha ,\kappa ,y} =
\int_{G_{\beta }} r(h)\delta _y(dh) =r(y)$$ for each continuous
bounded function $r$ on $G_{\beta }$. Therefore, $$S_{\alpha
}f_{\alpha } = p_{\alpha } \frac{\rho^l_{\mu _{\alpha }}
(\theta^{\beta }_{\alpha }(a_{\beta }),e_{\alpha })}{ \rho^l_{\mu
_{\alpha }} (\theta^{\beta }_{\alpha }(a_{\beta }),e_{\alpha
})}{\hat U}_{a_{\beta }} f_{\alpha }= p_{\alpha }{\hat U}_{a_{\beta
}} f_{\alpha }$$ on $L^1_{G_{\beta }}(G_{\alpha })$ for each $\alpha
$, where $\beta =\phi (\alpha )$. Moreover,
$$\| S \| = \sup_{\alpha \in \Lambda ; \beta =\phi (\alpha )}
p_{\alpha } \| {\hat U}_{a_{\beta }} \|_{L(X,X)}$$ with $X=L^{\infty
}(L^1_{G_{\beta }}(G_{\alpha }): \alpha <\beta \in \Lambda )$, where
$L(X,Y)$ denotes the normed space of all bounded linear operators
from $X$ to $Y$ with $X$ and $Y$ being complex normed spaces.
\par {\bf 14. Lemma.} {\it Let $\hat K$ be the scalar continuous operator
\par $(1)$ ${\hat K}f=pf$ for every $f\in L^{\infty }(L^1_{G_{\beta }}(G_{\alpha }): \alpha
<\beta \in \Lambda )$, that is $K_{\alpha }f_{\alpha }=p_{\alpha
}f_{\alpha }$ with $p_{\alpha }>0$  for each $\alpha \in \Lambda $.
Then this operator $\hat K$ satisfies Conditions 11$(1-3)$ and \par
$(2)$ $f\sim {\hat K}f$, that is $f\sim u$ by the definition means
that for every $t$ if $f_{\alpha }\ge 0$ and $u_{\alpha }\ge 0$ and
$t_{\alpha }\ge 0$ for each $\alpha \in \Lambda $ then $\{ f_{\alpha
}\wedge t_{\alpha }=0 \} \Leftrightarrow \{ u_{\alpha }\wedge
t_{\alpha }=0 \} $ for any $\alpha \in \Lambda $, where $f, u, t \in
L^{\infty }(L^1_{G_{\beta }}(G_{\alpha }): \alpha <\beta \in \Lambda
)$.}
\par {\bf Proof.} Properties 11$(1-3)$ are evidently satisfied for
$\hat K$. Condition $(2)$ is also fulfilled, since $supp (f_{\alpha
}) = supp (p_{\alpha }f_{\alpha })$ for each $\alpha \in \Lambda $.
\par {\bf 15. Theorem.} {\it Topological group rings
$L^{\infty }(L^1_{G_{\beta }}(G_{\alpha }, \mu_{\alpha }): \alpha
<\beta \in \Lambda )$ and $L^{\infty }(L^1_{G_{\beta }}(G_{\alpha },
\nu_{\alpha }): \alpha <\beta \in \Lambda )$ are isomorphic if and
only if measures $\mu _{\alpha }$ and $\nu _{\alpha }$ are
equivalent for each $\alpha \in \Lambda $.}
\par {\bf Proof.} If measures are equivalent, an isomorphism of topological group rings (and
algebras) is given by $$(1)\quad L^1_{G_{\beta }}(G_{\alpha },\mu
_{\alpha })\ni f_{\alpha }\mapsto f_{\alpha }\frac{d\mu _{\alpha
}}{d\nu _{\alpha }}\in L^1_{G_{\beta }}(G_{\alpha },\nu _{\alpha
})\quad \forall \alpha \in \Lambda .$$
\par Vice versa if topological group rings are isomorphic, all their
representations in $L(X,X)$ are equivalent, where $X$ is a complex
Banach space, $L(X,X)$ denotes the Banach space of all continuous
linear operators on $X$ into $X$. Particularly, ring representations
induced by unitary regular representations of groups $G_{\beta }$
are also equivalent. A regular unitary representation $T^{\beta ,\mu
_{\alpha }}: G_{\beta }\to U(X_{\alpha })$ is prescribed by the
formula
$$(2)\quad T^{\beta ,\mu _{\alpha }}(h) f_{\alpha }(x)=
\sqrt{\rho_{\mu _{\alpha }}(\theta ^{\beta }_{\alpha }(h),x)}
f_{\alpha }(\theta ^{\beta }_{\alpha }(h^{-1})x),$$ where $\beta
=\phi (\alpha )$, $~\alpha \in \Lambda $, $h\in G_{\beta }$,
$~X_{\alpha } := L^2(G_{\alpha },\mu_{\alpha },{\bf C})$,
$~U(X_{\alpha })$ denotes the unitary group on the Hilbert space
$X_{\alpha }$. The representation $T^{\beta ,\mu _{\alpha }}$ is
strongly continuous on each $G_{\beta }$ (see also
\cite{lujms150:4:08} and References 55 and 181 and 195 there and \S
10 above).
\par The family of all simple functions of the form
$$(3)\quad f_{\alpha } = \sum_{j=1}^n b_j \chi_{B_{\alpha ,j}}$$
with $n\in {\bf N}$, $~b_j\in {\bf C}$, open subsets $B_{\alpha ,j}$
in $(G_{\alpha },\tau _{\alpha })$, $B_{\alpha ,j}\cap B_{\alpha
,k}=\emptyset $ for all $j\ne k$, is dense in $X_{\alpha }$, since
the measure $\mu _{\alpha }$ is Borel regular. On the other hand,
$B_{\alpha, j}\cap \theta ^{\beta }_{\alpha }(G_{\beta })$ are open
in $G_{\beta }$, since $\tau _{\alpha }\cap {\theta ^{\beta
}_{\alpha }(G_{\beta })} \subset \theta ^{\beta }_{\alpha }(\tau
_{\beta })$. From the latter property it follows that
$$f_{\alpha }\circ \theta ^{\beta }_{\alpha } = \sum_{j=1}^n b_j \chi_{K_{\beta ,j}},$$
where $K_{\beta ,j}$ is open in $G_{\beta }$ for each $j$.
Therefore, the topological density $d(X_{\alpha })$ of the Hilbert
space $X_{\alpha }$ is not greater than that of $X_{\beta }$,
consequently, there exists an isometric linear embedding \par $(4)$
$\eta ^{\alpha }_{\beta }: X_{\alpha }\to X_{\beta }$.
\par Let $Y$ be a Banach space consisting of all vectors $y=(y_{\alpha }: ~
y_{\alpha }\in X_{\alpha }, ~ \alpha \in \Lambda )$ with $$\| y \| =
\sup_{\alpha \in \Lambda } \| y_{\alpha } \| _{X_{\alpha }}<\infty
,$$ it can be also denoted by $Y=l_{\infty }(X_{\alpha }: \alpha \in
\Lambda )$. The embedding $X_{\alpha }\oplus X_{\beta
}\hookrightarrow Y$ induces the embedding $L(X_{\alpha },X_{\alpha
})\oplus L(X_{\beta },X_{\beta })\hookrightarrow L(Y,Y)$. If a
strongly continuous unitary representation $T^{\alpha }: G_{\alpha
}\to L(X_{\alpha },X_{\alpha })$ is given for each $\alpha \in
\Lambda $, then taking Bochner's integral
$$R^{\alpha }(f_{\alpha }) y_{\alpha } := \int_{G_{\alpha }} f_{\alpha }
(g)T^{\alpha }(g)y_{\alpha } \mu _{\alpha }(dg),$$ where $f_{\alpha
}\in L^1_{G_{\beta }}(G_{\alpha })$ with $\beta =\phi (\alpha )$ and
$y_{\alpha }\in X_{\alpha }$, we get a strongly continuous
representation $R: L^{\infty }(L^1_{G_{\beta }}(G_{\alpha },
\mu_{\alpha }): \alpha <\beta \in \Lambda )\to L(Y,Y),$ since
$X_{\alpha }\oplus X_{\beta }$ has an embedding into $Y$ and
$$R(u_{\beta }{\tilde *}f_{\alpha }) (y_{\beta }\oplus y_{\alpha }) := R[\int_{G_{\beta }}
u_{\beta }(h)f_{\alpha }(\theta ^{\beta }_{\alpha }(h)g) \mu _{\beta
}(dh) ] y_{\beta }\oplus y_{\alpha } = $$ $$\int_{G_{\beta }}
\int_{G_{\alpha }} (T^{\beta }(h)u_{\beta }(h) y_{\beta })\oplus
(f_{\alpha } (\theta ^{\beta }_{\alpha }(h)g)T^{\alpha }(\theta
^{\beta }_{\alpha }(h)g)y_{\alpha }\mu_{\beta }(dh) \mu _{\alpha
}(dg)$$
$$=R^{\beta }(u_{\beta }){\tilde *} R^{\alpha }(f_{\alpha })
(y_{\beta }\oplus y_{\alpha })$$ for each $\alpha \in \Lambda $ with
$\beta = \phi (\alpha )$. Therefore, $R(u{\tilde \star }f)= R(u)
{\tilde \star }R(f)$ for each $f, u \in L^{\infty }(L^1_{G_{\beta
}}(G_{\alpha }, \mu_{\alpha }): \alpha <\beta \in \Lambda )$.
\par On the other hand, the embedding $\eta ^{\alpha }_{\beta }: X_{\alpha }\to X_{\beta }$
provides an equivalence relation $\Sigma ^{\alpha }$ so that
$T^{\beta ,\mu _{\alpha }}\times T^{\gamma ,\nu _{\beta }}$ induce a
unitary representation $T^{\beta , \gamma }: G_{\beta }\times
G_{\gamma }\to U(X_{\beta })$ for which $$(5)\quad T^{\beta , \gamma
}(g_{\beta } ,g_{\gamma }) = {\hat T}^{\beta ,\mu _{\alpha
}}(g_{\beta }) T^{\gamma ,\mu _{\beta }}(g_{\gamma })\in U(X_{\beta
}),$$ where a representation ${\hat T}^{\beta ,\mu _{\alpha }}$ is
induced by $T^{\beta ,\mu _{\alpha }} $ and the embedding $\eta
^{\alpha }_{\beta }$. Thus using an approximation of Dirac's measure
$\delta _z$ on $G_{\beta }$ we get that an equivalence of two
representations of two rings $L^{\infty }(L^1_{G_{\beta }}(G_{\alpha
}, \mu_{\alpha }): \alpha <\beta \in \Lambda )$ and $L^{\infty
}(L^1_{G_{\beta }}(G_{\alpha }, \nu_{\alpha }): \alpha <\beta \in
\Lambda )$ provides an intertwining operator $A^{\beta }$ of two
regular unitary representations $T^{\beta ,\mu _{\alpha }}$ and
$T^{\beta ,\nu _{\alpha }}$ such that $A^{\beta }: L^2(G_{\alpha
},\nu _{\alpha },{\bf C})\to L^2(G_{\alpha },\mu _{\alpha },{\bf
C})$ is a linear isomorphism of Hilbert spaces for each $\alpha \in
\Lambda $. Measures $\mu _{\alpha }$ and $\nu _{\alpha }$ are
regular Borel measures, while each group $G_{\alpha }$ is Hausdorff,
hence these regular representations distinguish different elements
of $G_{\beta }$, that is they are injective. Thus $$(6)\quad
T^{\beta ,\mu_{\alpha }} = A^{\beta }T^{\beta ,\nu _{\alpha
}}(A^{\beta })^{-1}.$$ From this it follows, as it is known from the
literature, that measures $\mu _{\alpha }$ and $\nu _{\alpha }$ are
equivalent (see also \cite{lujms150:4:08} and References 55 and 181
and 195 there). We shortly recall a way of the proof.
\par A Borel subset $J\in {\cal B}(G_{\alpha })$ is of $\mu _{\alpha }$
measure zero, i.e. $\int_{G_{\alpha }} \chi_J d\mu _{\alpha }=0$, if
and only if $\int_{G_{\alpha }} k\chi_J d\mu _{\alpha }=0$ for each
nonnegative continuous function $k$ on $G_{\alpha }$. If $V$ is a
linear topological isomorphism of $L^2(G_{\alpha },\mu_{\alpha
},{\bf C})$ onto itself, then $\int_{G_{\alpha }} V(\chi_J) d\mu
_{\alpha }=0$. Therefore, an operator $V$ preserves invariant a set
of nonnegative functions $s$ on $G_{\alpha }$ with $\mu_{\alpha }(s)
:= \int_{G_{\alpha }} s d\mu _{\alpha }=0$, that is the family of
all subsets in $G_{\alpha }$ of $\mu_{\alpha }$-measure zero is
invariant under $V$. \par Consider matrix elements
$$(8)\quad (T^{\beta ,\mu_{\alpha }}(g)x,y) = (A^{\beta }T^{\beta
,\nu _{\alpha }}(g)(A^{\beta })^{-1}x,y) $$ for each $x, y \in
X_{\alpha }$. Evidently, $\| y \| =0$ if and only if the scalar
product in $(8)$ is zero for each $x\in X_{\alpha }$ and $g\in
G_{\beta }$. Sets $E_{n,\mu _{\alpha }}$ and $E_{n,\nu _{\alpha }}$
for measures $\mu _{\alpha }$ and $\nu _{\alpha }$ can be considered
as $E_n$ in \S 10. Then the limit $$\lim_n \int_{G_{\alpha }}
|(A^{\beta })^{-1}[1-\chi _{E_{n,\mu _{\alpha }}}]| d\nu _{\alpha
}=0$$ exists, since $\lim_n \mu _{\alpha } (G_{\alpha }\setminus
E_{n,\mu _{\alpha }})=0$ and the linear operator $(A^{\beta })^{-1}$
is continuous. Symmetrically $$\lim_n \int_{G_{\alpha }} |A^{\beta
}[1-\chi _{E_{n,\nu _{\alpha }}}]| d\mu _{\alpha }=0.$$
\par Each left shift $L_h: G_{\alpha }\to G_{\alpha }$ with $h\in
\theta ^{\beta }_{\alpha }(G_{\beta })$, $L_hg=hg$, induces an
isometry $f_{\alpha }(g)\mapsto f_{\alpha }(hg)$ of $L^1_{G_{\beta
}}(G_{\alpha }, \mu _{\alpha })$ onto itself and also for
$L^1_{G_{\beta }}(G_{\alpha }, \nu _{\alpha })$. If $\{ U_{\alpha
,v}: v \} $ is a base of neighborhoods of the unit element in
$G_{\alpha}$, an arbitrary element $q\in G_{\alpha }$ is marked,
then there are elements $g_{\beta ,v,q}\in \theta ^{\beta }_{\alpha
}(G_{\beta })$ such that $g_{\beta ,v,q}U_{\alpha ,v}$ is a base of
neighborhoods of $q$.
\par By Cauchy-Bounyakovskii's inequality $|(T^{\beta ,\mu_{\alpha
}}(g)x_m,y)|\le \| x_m \| \| y \| $, since $\| T^{\beta ,\mu_{\alpha
}}(g) \| =1$, hence if $x_m\to 0$, then $(T^{\beta ,\mu_{\alpha
}}(g)x_m,y)\to 0$ uniformly by $g\in G_{\alpha }$ with $m$ tending
to the infinity. If $z_m \in L^2(G_{\alpha },\mu _{\alpha },{\bf
C})$ are nonzero vectors, then $b_mz_m/\| z_m\|\to 0$ for each
sequence of complex numbers $b_m$ tending to zero, when $m$ tends to
the infinity. A nonequivalence of measures would lead to a
contradiction when one regular representation would be strongly
continuous and another not on certain vectors, but these
representations are equivalent and related by Formula $(8)$. In view
of Lemma 12 and Formulas $(2,8)$ for each Borel subset $B$ in
$G_{\alpha }$: $ ~ \mu_{\alpha }(B)=0$ if and only if $\nu _{\alpha
}(B)=0$, consequently, measures $\mu _{\alpha }$ and $\nu _{\alpha
}$ are equivalent, since measures $\mu _{\alpha }$ and $\nu _{\alpha
}$ are Borel regular.
\par {\bf 16. Theorem.} {\it Let $G = \prod_{\alpha \in \Lambda }G_{\alpha }$ and
$H = \prod_{\alpha \in \Lambda }H_{\alpha }$ be two topological
groups supplied with box topologies $\tau ^b_G$ and $\tau ^b_H$
respectively, where topological groups $G_{\alpha }$ and $H_{\alpha
}$ for each $\alpha \in \Lambda $ satisfy Conditions 1$(1-4)$,
measures $\mu _{\alpha }$ on $G_{\alpha }$ and $\nu _{\alpha }$ on
$H_{\alpha }$ satisfy Conditions 2$(1-4)$, a directed set $\Lambda $
has not a minimal element. \par 1. If topological groups $G_{\alpha
}$ and $H_{\alpha }$ for each $\alpha \in \Lambda $ are
topologically isomorphic, then equivalent measures $\mu _{\alpha }$
and $\nu _{\alpha }$ exist so that topological algebras $L^{\infty
}(L^1_{G_{\beta }}(G_{\alpha },\mu_{\alpha }): \alpha <\beta \in
\Lambda )$ and $L^{\infty }(L^1_{H_{\beta }}(H_{\alpha },\nu
_{\alpha }): \alpha <\beta \in \Lambda )$ are isomorphic and their
isomorphism $\hat T$ satisfies properties $(1-3)$ below.
\par 2. If a bijective surjective
continuous mapping ${\hat T}$ of $L^{\infty }(L^1_{G_{\beta
}}(G_{\alpha }): \alpha <\beta \in \Lambda )$ onto $L^{\infty
}(L^1_{H_{\beta }}(H_{\alpha }): \alpha <\beta \in \Lambda )$ exists
and ${\hat T}^{-1}$ is continuous such that \par $(1)$ a mapping
${\hat T}= ({\hat T}_{\alpha }f_{\alpha }: \alpha \in \Lambda )$ is
linear so that ${\hat T}_{\alpha }: L^1_{G_{\beta }}(G_{\alpha })\to
L^1_{H_{\beta }}(H_{\alpha })$ for every $\alpha \in \Lambda $ with
$\beta = \phi (\alpha )$;
\par $(2)$ $\hat T$ is positive, that is $f_{\alpha }\ge 0$ in
$L^1_{G_{\beta }}(G_{\alpha })$ if and only if ${\hat T}_{\alpha
}f_{\alpha }\ge 0$ in $L^1_{H_{\beta }}(H_{\alpha })$;
\par $(3)$ ${\hat T}$ is a ring homomorphism, that is
${\hat T}(f{\tilde \star }u)= (f {\tilde \star } {\hat T}u)$ for
each $f, u \in L^{\infty }(L^1_{G_{\beta }}(G_{\alpha }): \alpha
<\beta \in \Lambda )$, \par then topological groups $G_{\alpha }$
and $H_{\alpha }$ are topologically isomorphic and measures $\mu
_{\alpha }$ and $\nu _{\alpha }$ are equivalent for each $\alpha \in
\Lambda $.}
\par {\bf Proof.} If $\omega _{\alpha }: G_{\alpha }\to H_{\alpha }$
is a topological group isomorphism for each $\alpha \in \Lambda $,
then an operator ${\hat T}$ with $({\hat T}_{\alpha }f_{\alpha })(h)
:= f_{\alpha }(\omega ^{-1}_{\alpha }(h))$ for every $\alpha \in
\Lambda $ and $h\in H_{\alpha }$ and $f_{\alpha }\in L^1_{G_{\beta
}}(G_{\alpha })$ has the desired properties. Taking a measure $\nu
_{\alpha }=\mu _{\alpha }\circ \omega _{\alpha }^{-1}$ on $H_{\alpha
}$ for any $\alpha \in \Lambda $ establishes an isometric
isomorphism ${\hat T}$ which satisfies Conditions $(1-3)$.
\par Conversely, let ${\hat T}$ satisfy the conditions of this theorem.
Then from the conditions of this theorem the algebraic isomorphisms
\par $(4)$ ${\hat U}(G)\cong {\hat U}(H)$ and ${\hat S}(G)\cong
{\hat S}(H)$ follow, where ${\hat U}(G)$ and ${\hat S}(G)$ denote
the families of all operators satisfying Conditions $11(1-3)$ and
(11$(1-3)$,14$(1,2)$) correspondingly. These algebraic homomorphisms
are induced by the operator $\hat T$ according to the formula
\par $(5)$ ${\hat U}(G)\ni {\hat P}\mapsto {\hat T}{\hat P}{\hat T}^{-1}\in
{\hat U}(H)$. \par In view of Theorem 11 and Lemma 14 there are
algebraic isomorphisms ${\hat G}\cong {\hat U}(G)/{\hat S}(G)$ and
${\hat H}\cong {\hat U}(H)/{\hat S}(H)$, since the group algebra
$L^{\infty }(L^1_{G_{\beta }}(G_{\alpha }): \alpha <\beta \in
\Lambda )$ is isomorphic with $L^{\infty }(L^1_{H_{\beta
}}(H_{\alpha }): \alpha <\beta \in \Lambda )$, where ${\hat G}$ and
${\hat H}$ denote the groups of left translation operators from \S
11. Therefore, algebraically the group ${\hat G}$ is isomorphic with
$G$ and the group ${\hat H}$ with $H$ respectively, since a directed
set $\Lambda $ has not a minimal element. But the mappings ${\hat
T}$ and ${\hat T}^{-1}$ are bijective, surjective and continuous,
applying $(3)$ and 11$(3,17,18,20)$ we get the topological
isomorphism $\omega _{\alpha }$ of $G_{\alpha }$ with $H_{\alpha }$
for each $\alpha \in \Lambda $. \par In view of Theorem 15 measures
$\mu_{\alpha }\circ \omega _{\alpha }^{-1}$ and $\nu _{\alpha }$ on
$H_{\alpha }$ are equivalent and the isometric isomorphism of group
algebras is provided by the mapping $$L^1_{G_{\beta }}(G_{\alpha
},\mu _{\alpha })\ni f_{\alpha }\mapsto (f_{\alpha }\circ
\omega_{\alpha }^{-1})\frac{d\mu _{\alpha }\circ \omega_{\alpha }
^{-1}}{d\nu _{\alpha }}\in L^1_{H_{\beta }}(H_{\alpha },\nu _{\alpha
})\quad \forall \alpha \in \Lambda .$$

\par {\bf 17. Definition.}
Let
$$(1)\quad (g_{\beta }*f_{\alpha })(s)=\int_{G_{\beta }} g_{\beta }(x)f_{\alpha
}(\theta ^{\beta }_{\alpha }(x^{-1})s)\mu_{\beta }(dx)\mbox{  and}$$
$$(2)\quad (g_{\beta }\check{*}f_{\alpha })(s)=\int_{G_{\beta }} g_{\beta }(x)f_{\alpha
}(s\theta ^{\beta }_{\alpha }(x^{-1}))\mu_{\beta }(dx)$$ for each
$s\in G_{\alpha }$.
\par The group algebra $L^{\infty }(L^1_{G_{\beta }}(G_{\alpha }): \alpha
<\beta \in \Lambda )$ will be called meta-commutative, if the
following condition is satisfied:
$$(3)\quad [\frac{f_{\beta }(x)}{\psi _{\beta
}(x^{-1}\theta ^{\gamma }_{\beta }(s))\rho^l_{\mu_{\beta }}(\theta
^{\gamma }_{\beta }(s^{-1}),\theta ^{\gamma }_{\beta }(s^{-1})x)}*
g_{\alpha }(\theta ^{\beta }_{\alpha }(x^{-1})\theta ^{\gamma
}_{\alpha }(s))](\theta ^{\gamma }_{\alpha }(s))$$ $$ = [g_{\beta
}(y)\check{*} f_{\alpha }(\theta ^{\gamma }_{\alpha }(s)\theta
^{\beta }_{\alpha }(y^{-1})](\theta ^{\gamma }_{\alpha }(s))$$ for
each $\alpha \in \Lambda $, every $s\in G_{\gamma }$ with $\beta
=\phi (\alpha )$, $ ~ \gamma =\phi (\beta )$, for every $f, g \in
L^{\infty }(L^1_{G_{\beta }}(G_{\alpha }): \alpha <\beta \in \Lambda
)$ such that $g_{\alpha }\circ \theta ^{\beta }_{\alpha }=g_{\beta
}$ and
\par $(4)$ $f_{\alpha }\circ \theta ^{\beta }_{\alpha }=f_{\beta }$
on $G_{\beta }$.
\par {\bf 18. Theorem.} {\it The group algebra $L^{\infty }(L^1_{G_{\beta }}(G_{\alpha }): \alpha
<\beta \in \Lambda )$ is meta-commutative if and only if a group $G$
is commutative.}
\par {\bf Proof.} A group $G$ is commutative if and only if a
group $G_{\alpha }$ is commutative for each $\alpha \in \Lambda $.
Since there are approximations of Dirac's measure $\delta _z$ on
$G_{\alpha }$ in this group algebra, then the group $G_{\alpha }$ is
commutative if and only if \par $(1)$ $f_{\alpha
}(y^{-1}s)=f_{\alpha }(sy^{-1})$ for each $f_{\alpha }\in
L^1_{G_{\beta }}(G_{\alpha })$ and any $x, s \in G_{\alpha }$. \par
On the other hand, each continuous bounded function $f_{\alpha }$ on
$G_{\alpha }$ is also continuous on $\theta ^{\beta }_{\alpha
}(G_{\beta })$ relative to the topology $\theta ^{\beta }_{\alpha
}(\tau _{\beta })$, since $\tau_{\alpha }\cap \theta ^{\beta
}_{\alpha }(G_{\beta })\subset \theta ^{\beta }_{\alpha }(\tau
_{\beta })$, consequently, $f_{\alpha}$ has a continuous bounded
restriction $f_{\beta }=f_{\alpha }\circ \theta^{\beta }_{\alpha } $
on the topological space $(G_{\beta },\tau _{\beta })$. This
restriction $f_{\beta }$ satisfies Condition 17$(4)$. The space of
bounded continuous functions $f_{\alpha }$ satisfying 17$(4)$ is
dense in $L^1_{G_{\beta }}(G_{\alpha })$. \par Therefore, it is
sufficient to demonstrate that Equality $(1)$ is equivalent to
17$(3)$ for any $y, s \in \theta^{\gamma }_{\alpha }(G_{\gamma })$,
since $G_{\gamma }$ is dense in $G_{\alpha }$ and in $G_{\beta }$,
where $\beta = \phi (\alpha )$ and $\gamma = \phi (\beta )$. From
Formulas 17$(1,2)$ we infer that
$$(2)\quad [\frac{f_{\beta }(x)}{\psi _{\beta
}(x^{-1}\theta ^{\gamma }_{\beta }(s))\rho^l_{\mu_{\beta }}(\theta
^{\gamma }_{\beta }(s^{-1}),\theta ^{\gamma }_{\beta }(s^{-1})x)}*
g_{\alpha }(\theta ^{\beta }_{\alpha }(x^{-1})\theta ^{\gamma
}_{\beta }(s))](\theta ^{\gamma }_{\beta }(s)) =$$
$$ \int_{G_{\beta }} \frac{f_{\beta }(x)}{\psi _{\beta
}(x^{-1}\theta ^{\gamma }_{\beta }(s))\rho^l_{\mu_{\beta }}(\theta
^{\gamma }_{\beta }(s^{-1}),\theta ^{\gamma }_{\beta
}(s^{-1})x)}g_{\alpha }(\theta ^{\beta }_{\alpha }(x^{-1})\theta
^{\gamma }_{\alpha }(s))\mu_{\beta }(dx)$$  $$= \int_{G_{\beta
}}g_{\beta }(y) f_{\alpha }(\theta ^{\gamma }_{\alpha }(s)\theta
^{\beta }_{\alpha }(y^{-1}))\frac{\psi _{\beta
}(y)\rho^l_{\mu_{\beta }}(\theta ^{\gamma }_{\beta
}(s^{-1}),y^{-1})}{\psi _{\beta }(y)\rho^l_{\mu_{\beta }}(\theta
^{\gamma }_{\beta }(s^{-1}),y^{-1})} \mu_{\beta }(dy)$$ after the
change of the integration variable $y=x^{-1}\theta ^{\gamma }_{\beta
}(s)$, where $x\in G_{\beta }$ and $s\in G_{\gamma }$. Therefore,
from Formulas $(2)$ and 17$(1,2)$ it follows, that Condition 17$(3)$
is equivalent to the equality
$$(3)\quad \int_{G_{\beta }} g_{\beta }(x)f_{\alpha
}(\theta ^{\beta }_{\alpha }(x^{-1})\theta ^{\gamma }_{\beta
}(s))\mu_{\beta }(dx) =\int_{G_{\beta }} g_{\beta }(x)f_{\alpha
}(\theta ^{\gamma }_{\beta }(s))\theta ^{\beta }_{\alpha
}(x^{-1}))\mu_{\beta }(dx).$$ A bounded continuous function
$g_{\beta }$ is arbitrary in $L^1_{G_{\gamma }}(G_{\beta })$ in the
latter formula, consequently, Equality $(1)$ is satisfied if and
only if 17$(3)$ is valid.

\end{document}